\documentclass[10 pt]{amsart}
\usepackage{amsfonts}
\usepackage{ifthen}
\usepackage{amsthm}
\usepackage{amsmath}
\usepackage{graphicx}
\usepackage{amscd,amssymb,amsthm}
\usepackage{graphicx}
\usepackage{epstopdf}
\usepackage{hyperref}

\newcounter{minutes}
\divide\time by 60
\newcounter{hours}
\multiply\time by 60 \addtocounter{minutes}{-\time}

\newtheorem{lemma}{Lemma}
\newtheorem{theorem}{Theorem}
\newtheorem{corollary}{Corollary}

\newtheorem{problem}{Open Problem}

\newcommand{\real}{\operatorname{Re}}
\newcommand{\imag}{\operatorname{Im}}

\keywords{Coulomb wave function; univalent, starlike functions; radius of univalence, starlikeness and convexity;
zeros of regular Coulomb wave functions; Mittag-Leffler expansions; Laguerre-P\'olya class of entire functions.}
\subjclass[2010]{30C45, 30C15, 33C10}

\begin{document}

\title{Radii of starlikeness and convexity of regular Coulomb wave functions}

\author[\'A. Baricz]{\'Arp\'ad Baricz$^{\bigstar}$}
\address{Department of Economics, Babe\c{s}-Bolyai University, Cluj-Napoca
400591, Romania}
\address{Institute of Applied Mathematics, \'Obuda University, 1034
Budapest, Hungary}
\email{bariczocsi@yahoo.com}

\author[M. \c{C}a\u{g}lar]{Murat \c{C}a\u{g}lar}
\address{Department of Mathematics, Faculty of Science and Letters, Kafkas
University, Kars 36100, Turkey} \email{mcaglar25@gmail.com}

\author[E. Deniz]{Erhan Deniz}
\address{Department of Mathematics, Faculty of Science and Letters, Kafkas
University, Kars 36100, Turkey} \email{edeniz36@gmail.com}

\author[E. Toklu]{Evrim Toklu}
\address{Department of Mathematics, Faculty of Science, Atat\"urk University,
25240 Erzurum, Turkey} \email{evrimtoklu@gmail.com}

\thanks{$^{\bigstar}$The research of \'A. Baricz was supported by a research grant of the Babe\c{s}-Bolyai University for young researchers with project number GTC-31777.}

\def\thefootnote{}
\footnotetext{ \texttt{File:~\jobname .tex,
          printed: \number\year-0\number\month-\number\day,
          \thehours.\ifnum\theminutes<10{0}\fi\theminutes}
} \makeatletter\def\thefootnote{\@arabic\c@footnote}\makeatother

\maketitle

\begin{center}
{Dedicated to the memory of Petru T. Mocanu}
\end{center}

\begin{abstract}
In this paper our aim is to determine the radii of univalence, starlikeness and convexity of the
normalized regular Coulomb wave functions for two different kinds of normalization. The key tools in
the proof of our main results are the Mittag-Leffler expansion for regular Coulomb wave functions, and properties of zeros
of the regular Coulomb wave functions and their derivatives. Moreover, by using the technique of differential
subordinations we present some conditions on the parameters of the regular Coulomb wave function in order
to have a starlike normalized form. In addition, by using the Euler-Rayleigh inequalities
we obtain some tight bounds for the radii of starlikeness of the normalized regular Coulomb wave functions.
Some open problems for the zeros of the regular Coulomb wave functions are also stated which may be of interest for further research.
\end{abstract}

\section{\bf Introduction}

Motivated by the beautiful properties of special functions, in the recent years there was a vivid interest on geometric properties of special functions, like Bessel, Struve and Lommel functions of the first kind; see the papers \cite{aktas,Ba4,bcd,bdy,bdoy,Ba0,Ba2,Ba1,Ba3,basz2,L1,szasz,Sza} and the references therein. The determination of the radii of starlikeness and convexity of some
normalized forms of these special functions was studied in details by the first author and his collaborators.
More precisely, recently the authors in \cite{Ba0}, \cite{Ba2} and \cite{Ba1} obtained the radii of starlikeness of order $\beta,$ the radii of convexity
of order $\beta $ and the radii of $\alpha$-convexity of order $\beta $ for three kind of normalized forms of the Bessel functions of the first kind
$J_{\nu}$ in the case when $\nu>-1.$ Moreover, recently the authors in \cite{Ba3} and \cite{Sza} investigated the radii
of starlikeness of order $\beta $ and the radii of convexity of order $\beta $ for the normalized
Bessel functions in the case when $\nu\in (-2,-1)$ by using a little bit more complicated analysis than in the
case when $\nu>-1.$ In these papers the positive zeros of the Bessel functions of the first kind played an important role, and the knowledge of the distribution
of the real and nonreal zeros of Bessel functions of the first kind was essential. It was also shown that the technique used in order to deduce the radii of starlikeness and convexity for the case of Bessel functions it is also working for normalized Struve and Lommel functions, but some basic facts for the positive zeros of these functions were needed. See the papers \cite{bdoy}, \cite{basz2} and \cite{L1} for more details. Here the Laguerre-P\'olya class of real entire functions played an important role. These radii of starlikeness and convexity are in fact solutions of some transcendental equations and hence it is also worth
to find some estimates for them. In a recent paper \cite{aktas} the authors found some tight lower and upper bounds for the radii of starlikeness of some normalized Bessel, Struve and Lommel functions and they also have shown that the corresponding radii of starlikeness are in fact equal to the radii of univalence. Motivated by the above series of papers in this subject, in this paper our aim is to present some similar results for the normalized forms of the regular Coulomb wave function. We are mainly interested in this paper on the determination of the radii of univalence, starlikeness and convexity; we will present also some conditions on the parameters such that the corresponding regular Coulomb wave functions belong to some special class of analytic functions; and finally, our aim is also to show some lower and upper bounds for the radii of univalence and starlikeness. As we can see the technique used for Bessel functions of the first kind is also working here, and thus the regular Coulomb wave functions belong also to the family of entire functions which have nice geometric properties. The paper is organized as follows: the rest of this section is dedicated to some basic definitions. Section 2 contains the results concerning the radii of univalence, starlikeness and convexity of normalized regular Coulomb wave functions, and its last part is devoted to the study of starlikeness of regular Coulomb wave functions. At the end of this section some open problems are stated for zeros of regular Coulomb wave functions, which may be of interest for further research. In section 3 we present the proofs of the preliminary and main results.

Before to start the presentation of the results we recall first some basic definitions. For $r>0$ we denote by $\mathbb{D}_r=\left\{z\in\mathbb{C}: |z|<r\right\}$ the open disk of radius $r$
centered at the origin. Let $f:\mathbb{D}_r\to\mathbb{C}$ be the function defined by%
\begin{equation}
f(z)=z+\sum_{n\geq 2}a_{n}z^{n},  \label{eq0}
\end{equation}%
where $r$ is less or equal than the radius of convergence of the above power series. We say that the
function $f,$ defined by (\ref{eq0}), is starlike in the disk $\mathbb{D}_r$
if $f$ is univalent in $\mathbb{D}_r$, and $f(\mathbb{D}_r)$ is a
starlike domain in $\mathbb{C}$ with respect to the origin. Analytically, the function $f$ is starlike in
$\mathbb{D}_r$ if and only if $$\real\left[\frac{zf^{\prime }(z)}{f(z)}%
\right]>0 \quad \mbox{for all}\ \ z\in
\mathbb{D}_r.$$ For $\beta \in [0,1)$ we say that
the function $f$ is starlike of order $\beta $ in $\mathbb{D}_r$ if and
only if $$\real\left[\frac{zf^{\prime }(z)}{f(z)}\right]>\beta \quad \mbox{for all}\ \ z\in
\mathbb{D}_r.$$
The real number
\begin{equation*}
r_{\beta }^{\ast }(f)=\sup \left\{ r>0 \left|\real%
\left[\frac{zf^{\prime }(z)}{f(z)}\right] >\beta \;\text{for all }z\in
\mathbb{D}_r\right.\right\}
\end{equation*}%
is called the radius of starlikeness of order $\beta $ of the function $f$. Note that $%
r^{\ast }(f)=r_{0}^{\ast }(f)$ is in fact the largest radius such that the image
region $f(\mathbb{D}_{r^{\ast}(f)})$ is a starlike domain with respect to
the origin.

The function $f,$ defined by (\ref{eq0}), is convex in the disk $\mathbb{D}_r$ if $f$
is univalent in $\mathbb{D}_r$, and $f(\mathbb{D}_r)$ is a
convex domain in $\mathbb{C}.$ Analytically, the function $f$ is convex in $\mathbb{D}_r$ if and only
if $$\real\left[ 1+\frac{zf^{\prime \prime }(z)}{f^{\prime }(z)}\right]
>0  \quad \mbox{for all}\ \ z\in
\mathbb{D}_r.$$
For $\beta \in[0,1)$ we say that the
function $f$ is convex of order $\beta $ in $\mathbb{D}_r$ if and only if $$
\real\left[1+\frac{zf^{\prime \prime }(z)}{f^{\prime }(z)}\right]
>\beta \quad \mbox{for all}\ \ z\in
\mathbb{D}_r.$$ The radius of convexity of order $\beta $
of the function $f$ is defined by the real number%
\begin{equation*}
r_{\beta }^{c}(f)=\sup \left\{ r>0 \left|\real\left[1+
\frac{zf^{\prime \prime }(z)}{f^{\prime }(z)}\right]>\beta \;\text{for all }%
z\in\mathbb{D}_r\right.\right\} .
\end{equation*}
Note that $r^{c}(f)=r_{0}^{c}(f)$ is the largest radius such that the image
region $f(\mathbb{D}_{r^{c}(f)})$ is a convex domain. Finally, we recall that the radius of univalence of the
analytic function $f$ in the form of \eqref{eq0} is the largest radius $r$ such that $f$ maps $\mathbb{D}_r$ univalently into $f(\mathbb{D}_r).$

\section{\bf The radii of univalence, starlikeness and convexity of Coulomb wave functions}
\setcounter{equation}{0}

Let $_{1}F_{1}$ denote the Kummer confluent hypergeometric function. The regular
Coulomb wave function is defined as
$$F_{L,\eta}(z)=z^{L+1}e^{-iz}C_{L}(\eta){}_{1}F_{1}(L+1-i\eta,2L+2;2iz)=C_{L}(\eta )\sum_{n\geq 0}a_{L,n}z^{n+L+1},$$
where $L,\eta \in \mathbb{C},$ $z\in\mathbb{C},$
\begin{equation*}
C_{L}(\eta )=\frac{2^{L}e^{-\frac{\pi \eta }{2}}\left\vert \Gamma (L+1+i\eta
)\right\vert }{\Gamma (2L+2)}=\left\{
\begin{array}{cc}
\frac{2^{L}}{(2L+1)!}\sqrt{\frac{2\pi \prod\nolimits_{k=0}^{L}(k^{2}+\eta
^{2})}{\eta (e^{2\pi \eta }-1)},} & if\text{ }\eta \neq 0 \\
\frac{2^{L}L!}{(2L+1)!}, & if\text{ }\eta =0%
\end{array}%
\right.
\end{equation*}%
and%
\begin{equation*}
a_{L,0}=1,\;\text{\ \ }a_{L,1}=\frac{\eta }{L+1},\text{ \ \ }a_{L,n}=\frac{%
2\eta a_{L,n-1}-a_{L,n-2}}{n(n+2L+1)},\text{ \ \ }n\in \{2,3,\dots\}.
\end{equation*}

In this paper we focus on the following normalized forms%
$$f_{L,\eta }(z)=\left[ C_{L}^{-1}(\eta )F_{L,\eta }(z)\right] ^{\frac{1}{L+1}}\quad \mbox{and}
\quad g_{L,\eta }(z)=C_{L}^{-1}(\eta )z^{-L}F_{L,\eta }(z),$$
and we are going to present some interesting results on these two functions by using some
basic techniques of the geometric function theory.

\subsection{The radii of starlikeness of the functions $f_{L,\eta}$ and $g_{L,\eta}$}
Our first main result is related to the radii of starlikeness of these two normalized regular
Coulomb wave functions. However, before to present our first main result we focus first
on a technical result, which will be needed in the sequel.

\begin{lemma}
\label{lem}
If $\lambda\in[0,1],$ $a>b>0$ and $z\in\mathbb{C}$ such that $|z|<b,$ then
\begin{equation}\label{ineqlem}
\lambda \real \frac{z^2}{a(a\pm z)} -\real \frac{z^2}{b(b\pm z)}\geq \lambda \frac{|z|^2}{a(a\pm|z|)}-\frac{|z|^2}{b(b\pm|z|)}.
\end{equation}
\end{lemma}

Now, we are going to determine the radii of starlikeness of the normalized Coulomb wave functions.

\begin{theorem}\label{t2} Let $L>-1,$ $\eta \leq 0$ and $\beta \in [0,1).$ Then the radius of starlikeness of order
$\beta $ of the function $f_{L,\eta}$ is the smallest positive root of the equation
$$rF_{L,\eta }^{\prime }(r)-\beta (L+1)F_{L,\eta }(r)=0.$$ Moreover, the radius of starlikeness of order $\beta $ of the function $g_{L,\eta
}$ is the smallest positive root of the equation $$(L+\beta )F_{L,\eta }(r)-rF_{L,\eta }^{\prime }(r)=0.$$
\end{theorem}

It is important to mention that the regular Coulomb wave function is actually a generalization of a transformation of the Bessel function of the first kind. Namely, we have the relation $$F_{L,0}(z)=\sqrt{\frac{\pi z}{2}}J_{L+\frac{1}{2}}(z),$$
where $J_L$ stands for the Bessel function of the first kind and order $L.$ Taking into account this, it is clear that Theorem \ref{t2} in particular when $\eta=0$ reduces to the following interesting results, and one of them naturally complements the results from \cite{Ba0}. Note that in the simplification of the expressions $f_{L,0}(z)$ and $g_{L,0}(z)$ we used the Legendre duplication formula for the Euler gamma function, that is, we used
$$C_L^{-1}(0)=\frac{\Gamma(2L+2)}{2^L\Gamma(L+1)}=\frac{1}{\sqrt{\pi}}2^{L+1}\Gamma\left(L+\frac{3}{2}\right).$$ The result on $f_{L-\frac{1}{2},0}$ is new and complements \cite[Theorem 1]{Ba0}, however, the result on $g_{L-\frac{1}{2},0}$ is not new, it was proved in \cite[Theorem 1]{Ba0}. Thus the second part of Theorem \ref{t2} is actually a generalization of part {\bf b} of \cite[Theorem 1]{Ba0}.

\begin{corollary}
If $L>-\frac{1}{2}$ and $\beta\in[0,1),$ then the radius of starlikeness of order $\beta$ of the function
$$z\mapsto f_{L-\frac{1}{2},0}(z)=\left[2^{L}\Gamma(L+1)\sqrt{z}J_L(z)\right]^{\frac{1}{L+\frac{1}{2}}}$$
is the smallest positive root of the equation $$rJ_L'(r)-\left[\beta\left(L+\frac{1}{2}\right)-\frac{1}{2}\right]J_L(r)=0.$$
Moreover, under the same conditions the radius of starlikeness of order $\beta$ of the function
$$z\mapsto g_{L-\frac{1}{2},0}(z)=2^{L}\Gamma(L+1)z^{1-L}J_L(z)$$
is the smallest positive root of the equation $$rJ_L'(r)-(\beta+L-1)J_L(r)=0.$$
\end{corollary}

\subsection{The radius of univalence of normalized Coulomb wave functions} In this subsection our aim is to show that the radii of univalence of the
functions $f_{L,\eta}$ and $g_{L,\eta}$ correspond to the radii of starlikeness. Moreover, we will present some lower and upper bounds for the radii of starlikeness of these functions. The method which we use goes back to Euler and it was used recently in \cite{aktas}, see also \cite{ismail,wilf} for more details.

\begin{theorem}\label{thuniv}
If $L>-1$ and $\eta\leq0,$ then the radii of univalence of the
functions $f_{L,\eta}$ and $g_{L,\eta}$ correspond to the radii of starlikeness. Moreover, if $L>-1$ and $\eta<0,$ then we have that
\begin{equation}\label{ray1}
r^{\ast}(f_{L,\eta})>\frac{(L+1)^2\sqrt{2L+3}}{\sqrt{L^4+6L^3+(\eta^2+12)L^2+2(3\eta^2+5)L+3(2\eta^2+1)}}
\end{equation}
and
\begin{equation}\label{ray2}
r^{\ast}(f_{L,\eta})<\frac{2(L + 2)(L + 1)^2(L^4 + 6L^3 + (\eta^2+12)L^2+2(3\eta^2+5)L+3(2\eta^2+1))}
{\eta(4L^7+\eta_1L^6+\eta_2L^5+8\eta_3L^4+2\eta_4L^3+\eta_5 L^2+9(5-28\eta^2)L-72\eta^2)},
\end{equation}
where $\eta_1=41-8\eta^2,$ $\eta_2=163-74\eta^2,$ $\eta_3=41-32\eta^2,$ $\eta_4=178-223\eta^2,$ $\eta_5=199-438\eta^2.$ In addition, the radius of starlikeness (and univalence) of $g_{L,\eta}$ satisfies
\begin{equation}\label{ray3}
r^{\ast}(g_{L,\eta})>\frac{(L+1)\sqrt{2L+3}}{\sqrt{3L^2+2(\eta^2+3)L+3(2\eta^2+1)}}
\end{equation}
and
\begin{equation}\label{ray4}
r^{\ast}(g_{L,\eta})<\frac{(L+1)(L+2)(3L^2+2(\eta^2+3)L+3(2\eta^2+1))}{\eta(8L^4+(31-16\eta^2)L^3+2(19-26\eta^2)L^2+3(5-22\eta^2)L-36\eta^2)}.
\end{equation}
\end{theorem}

\subsection{The radius of convexity of the function $f_{L,\eta}$ and $g_{L,\eta}$}

The next set of results concerns the radii of convexity of the functions $f_{L,\eta}$ and $g_{L,\eta}.$

\begin{theorem}\label{thconv}
If $L>-\frac{1}{2}$ and $\eta\leq0,$ then the radius of convexity of order $\beta\in[0,1)$ of $f_{L,\eta}$ is the smallest positive root of the equation
\begin{equation}\label{eqconv}1+\frac{%
rF_{L,\eta }^{\prime \prime }(r)}{F_{L,\eta }^{\prime }(r)}-\frac{L}{L+1}\frac{rF_{L,\eta }^{\prime }(r)}{F_{L,\eta }(r)}=\beta.\end{equation}
Moreover, if $L>-1$ and $\eta\leq 0,$ then the radius of convexity of order $\beta$ of $g_{L,\eta}$ is the smallest positive root of the equation
$$r^2F_{L,\eta}''(r)-(2L+\beta-1)rF_{L,\eta}'(r)+L(L+\beta)F_{L,\eta}(r)=0.$$
\end{theorem}

It is worth to mention that when $\eta=0$ the above result on $f_{L,\eta}$ reduces to the following new result on Bessel functions, which naturally complement the results from \cite{Ba1}. The result on $g_{L,\eta}$ in the particular case $\eta=0$ reduce to a known result from \cite[Theorem 1.2]{Ba1}.

\begin{corollary}
If $L>0$ and $\beta\in[0,1),$ then the radius of convexity of order $\beta$ of the function $f_{L-\frac{1}{2},0}$ is the smallest positive root of the transcendental equation
$$1+\frac{r^2J_L''(r)+rJ_L'(r)-\frac{1}{4}J_L(r)}{rJ_L'(r)+\frac{1}{2}J_L(r)}-\frac{L}{L+1}\frac{rJ_L'(r)+\frac{1}{2}J_L(r)}{J_L(r)}=\beta.$$
Moreover, if $L>-\frac{1}{2}$ and $\beta\in[0,1),$ then the radius of convexity of order $\beta$ of the function $g_{L-\frac{1}{2},0}$ is the smallest positive root of the equation
$$r^2J_L''(r)+(3-2L-\beta)rJ_L'(r)+(L-1)(L+\beta-1)J_L(r)=0.$$
\end{corollary}

\subsection{Starlikeness of the regular Coulomb wave function}

In this subsection our aim is to present a general result on the starlikeness of the normalized regular Coulomb wave function $g_{L,\eta}.$ In order to do this we shall use the technique of differential subordinations by Miller and Mocanu, see \cite{Mil1} and \cite{Mil2} for more details. In what follows we use the notation $\mathbb{D}=\mathbb{D}_1=\left\{\left.z\in\mathbb{C}\right||z|<1\right\}.$

\begin{lemma}\label{l1}
Let $\Omega\subseteq\mathbb{C}$ be a set in the complex plane $\mathbb{C}$ and $\Psi:
\mathbb{C}^{3}\times \mathbb{D}\to\mathbb{C}$ a function, that satisfies the admissibility condition $\Psi \left( \rho
i,\sigma ,\mu +i\nu;z\right) \notin \Omega,$ where $z\in \mathbb{D}$, $\rho ,\sigma ,\mu
,\nu\in\mathbb{R}$ with $\mu +\sigma \leq 0$ and $\sigma \leq -(1+\rho ^{2})/2.$ If $h$ is
analytic in the unit disk $\mathbb{D},$ with $h(0)=1$ and $\Psi \left(
h(z),zh^{\prime }(z),z^{2}h^{\prime \prime }(z);z\right) \in \Omega$ for all $%
z\in \mathbb{D},$ then $\real h(z)>0$ for all $z\in \mathbb{D}.$ In particular, if we only
have $\Psi:\mathbb{C}^{2}\times \mathbb{D}\to\mathbb{C},$ the admissibility condition reduces
to $\Psi \left( \rho i,\sigma
;z\right) \notin \Omega$ for all $z\in \mathbb{D}$ and $\rho ,\sigma \in\mathbb{R}$ with $\sigma \leq -(1+\rho ^{2})/2.$
\end{lemma}

Now, we are ready to state the main result of this subsection.

\begin{theorem}\label{thsub}
If $\eta ,L\in\mathbb{C}$ are such that $\real L\geq \frac{1}{2},$ $\imag L\geq 1$ and $\left(1+\imag L+\left\vert \eta \right\vert \right)
^{2}\leq \left(\real L-\frac{1}{2}\right) ^{2}$, then $\real g_{L,\eta}(z)>0$ for all $z\in\mathbb{D}.$ Moreover, if $\left\vert \eta \right\vert \leq \real L-\frac{1}{3}\left(\imag L\right) ^{2}-\frac{1}{4}$, then $g_{L,\eta}$ map the open unit disk $\mathbb{D}$ into a starlike domain with respect to origin.
\end{theorem}

\subsection{Some open problems related to regular Coulomb wave functions} In the process of writing this paper we faced some problems which we were not able to solve. A solution of these problems would enable us to complete the picture on the geometric properties of the regular Coulomb wave functions. The first two open problems are related to the $n$th zero $\rho_{L,\eta,n}$ of the regular Coulomb wave function $F_{L,\eta}$ and are motivated by the fact that the zeros of the Bessel function of the first kind are increasing with respect to the order. A positive answer to these open problems would help us to find sufficient and necessary conditions on the parameters $L$ and $\eta$ such that the normalized forms $f_{L,\eta}$ and $g_{L,\eta}$ belong to the class of univalent, starlike or convex functions. Moreover, it would be possible to find also sufficient and necessary conditions such that $g_{L,\eta}$ have close-to-convex derivatives.

\begin{problem}
If $\eta\in\mathbb{R}$ and $n\in\mathbb{N}$ are fixed, is it true that $L\mapsto \rho_{L,\eta,n}$ is increasing on $(-1,\infty)$?
\end{problem}

In the following open problem the word ''somehow'' means that we would like to see the derivative of the zeros with respect to the order in terms of the regular Coulomb wave function values, or its zeros, or in terms of integrals of which kernel contains some known special function.

\begin{problem}
Is it possible to express somehow the derivative $\partial \rho_{L,\eta,n}/\partial L$?
\end{problem}

The next open problem is related to the radius of convexity of the function $g_{L,\eta}.$ In the proof of Theorem \ref{thconv} we used the fact that
 for $L>-1$ and $\eta\in\mathbb{R}$ the function $r\mapsto rF_{L,\eta}'(r)-LF_{L,\eta}(r)$ has only real zeros. Now we are asking what happens in general with the zeros of the linear combinations of the normalized Coulomb wave function and its derivative.

\begin{problem}
Under which conditions on $L,$ $\eta$ and $\alpha$ the function $r\mapsto rF_{L,\eta}'(r)+\alpha F_{L,\eta}(r)$ has only real zeros?
\end{problem}

In the proof of the main results concerning the radii of univalence, starlikeness and convexity it was essential the fact that $\eta\leq0.$ The last open problem of this paper is related to this condition.

\begin{problem}
Is it possible to relax the condition $\eta\leq 0$ in the main results?
\end{problem}

\section{\bf Proofs of the main results}
\setcounter{equation}{0}

\begin{proof}[\bf Proof of Lemma \ref{lem}]
Observe that it is enough to show the inequality for the sign $-$ or $+$ in the brackets since one implies the other by changing $z$ to $-z.$ Taking into account this we shall prove the inequality
\begin{equation}\label{ineqlem1}\lambda \real \frac{z^2}{a(a-z)} -\real \frac{z^2}{b(b-z)}\geq \lambda \frac{|z|^2}{a(a-|z|)}-\frac{|z|^2}{b(b-|z|)},\end{equation}
where $\lambda\in[0,1],$ $a>b>0$ and $z\in\mathbb{C}$ such that $|z|<b.$ For $z=x+iy$ and $m=|z|$ we consider the function $\varphi:[m,\infty)\to\mathbb{R},$ defined by
$$\varphi(t)=\real \frac{z^2}{t(t-z)}-\frac{|z|^2}{t(t-|z|)}=\frac{m-x}{t}-\frac{t}{t-m}+\frac{t^2-tx}{t^2-2tx+m^2}.$$
We will show that this function maps the interval $[m,\infty)$ into $(-\infty,0]$ and increases. For $t>m$ we have
\begin{align*}
\varphi'(t)&=-\frac{m-x}{t^2}+\frac{m}{(t-m)^2}+\frac{2tm^2-m^2x-t^2x}{(t^2-2tx+m^2)^2}\\
&\geq -\frac{m-x}{(t-m)^2}+\frac{m}{(t-m)^2}+\frac{2tm^2-m^2x-t^2x}{(t^2-2tx+m^2)^2}\\
&= \frac{x}{(t-m)^2}+\frac{2tm^2-m^2x-t^2x}{(t^2-2tx+m^2)^2}\\
&\geq \frac{x}{t^2-2tx+m^2}+\frac{2tm^2-m^2x-t^2x}{(t^2-2tx+m^2)^2}\\
&=\frac{2t(m^2-x^2)}{(t^2-2tx+m^2)^2}>0,
\end{align*}
which implies that $\varphi$ is indeed increasing. Moreover, since $\lim_{t\searrow m}\varphi(t)=-\infty$ and $\lim_{t\to\infty}\varphi(t)=0,$ it follows that
$\varphi$ maps indeed the interval $[m,\infty)$ into $(-\infty,0].$ In other words, for $|z|<\mu$ we have
$$\real \frac{z^2}{\mu(\mu-z)}\leq\frac{|z|^2}{\mu(\mu-|z|)},$$
which implies that
$$\lambda\left[\real\frac{z^2}{a(a-z)}-\frac{|z|^2}{a(a-|z|)}\right]\geq \real\frac{z^2}{a(a-z)}-\frac{|z|^2}{a(a-|z|)},$$
where $|z|<a.$ On the other hand, since $\varphi$ is increasing we have that if $a> b> m,$ then $\varphi(a)> \varphi(b),$ that is,
$$\real\frac{z^2}{a(a-z)}-\frac{|z|^2}{a(a-|z|)}> \real\frac{z^2}{b(b-z)}-\frac{|z|^2}{b(b-|z|)}.$$ Combining these two inequalities we have that
\eqref{ineqlem1} is true, and this completes the proof.
\end{proof}

\begin{proof}[\bf Proof of Theorem \ref{t2}]
The Weierstrassian canonical product expansion of the Coulomb wave function reads as (see \cite{St})
\begin{equation}\label{prodco}F_{L,\eta }(z)=C_{L}(\eta )z^{L+1}e^{\frac{\eta z}{L+1}}\prod\limits_{n\geq
1}\left( 1-\frac{z}{\rho _{L,\eta ,n}}\right) e^{\frac{z}{\rho _{L,\eta ,n}}},\end{equation}
where $\rho _{L,\eta ,n}$ is the $n$th zero of the Coulomb wave function. Logarithmic
derivation yields
$$\frac{F_{L,\eta }^{\prime }(z)}{F_{L,\eta }(z)}=\frac{L+1}{z}+\frac{\eta }{%
L+1}-\sum\limits_{n\geq 1}\frac{z}{\rho _{L,\eta ,n}(\rho _{L,\eta ,n}-z)},$$
and this implies that
$$\frac{zf_{L,\eta }^{\prime }(z)}{f_{L,\eta }(z)}=\frac{1}{L+1}\frac{zF_{L,\eta}'(z)}{F_{L,\eta}(z)}=1+\frac{\eta z}{\left(
L+1\right) ^{2}}-\frac{1}{L+1}\sum\limits_{n\geq 1}\frac{z^{2}}{\rho
_{L,\eta ,n}(\rho _{L,\eta ,n}-z)}$$
and
$$\frac{zg_{L,\eta }^{\prime }(z)}{g_{L,\eta }(z)}=-L+\frac{zF_{L,\eta}'(z)}{F_{L,\eta}(z)}=1+\frac{\eta z}{L+1}%
-\sum\limits_{n\geq 1}\frac{z^{2}}{\rho _{L,\eta ,n}(\rho _{L,\eta ,n}-z)}.$$
Now, taking into account that the zeros $\rho _{L,\eta ,n}$ can be separated into
positive $x_{L,\eta ,n}$ and negative $y_{L,\eta ,n}$ zeros the above equations can be written as (for the first equation see \cite{Ba4})
$$\frac{F_{L,\eta }^{\prime }(z)}{F_{L,\eta }(z)}=\frac{L+1}{z}+\frac{\eta }{%
L+1}-\sum\limits_{n\geq 1}\left[\frac{z}{x_{L,\eta ,n}(x_{L,\eta ,n}-z)}+\frac{z}{y_{L,\eta ,n}(y_{L,\eta ,n}-z)}\right],$$
$$\frac{zf_{L,\eta }^{\prime }(z)}{f_{L,\eta }(z)}=1+\frac{\eta z}{\left(
L+1\right) ^{2}}-\frac{1}{L+1}\sum\limits_{n\geq 1}\left[\frac{z^{2}}{x_{L,\eta ,n}(x_{L,\eta ,n}-z)}+\frac{z^{2}}{y_{L,\eta ,n}(y_{L,\eta ,n}-z)}\right]$$
and
$$\frac{zg_{L,\eta }^{\prime }(z)}{g_{L,\eta }(z)}=1+\frac{\eta z}{L+1}%
-\sum\limits_{n\geq 1}\left[\frac{z^{2}}{x_{L,\eta ,n}(x_{L,\eta ,n}-z)}+\frac{z^{2}}{y_{L,\eta ,n}(y_{L,\eta ,n}-z)}\right].$$
In order to prove the assertions of the theorem we need to show that for the corresponding values of $L,\eta $ and $\beta $
the inequalities
\begin{equation}
\real\left[\frac{zf_{L,\eta }^{\prime }(z)}{f_{L,\eta }(z)}\right]
>\beta \text{ \ \ and \ \ }\real\left[\frac{zg_{L,\eta }^{\prime }(z)}{%
g_{L,\eta }(z)}\right]>\beta  \label{eq6}
\end{equation}%
are valid for $z\in \mathbb{D}_{r_{\beta }^{\ast }(f_{L,\eta })}$ and $z\in
\mathbb{D}_{r_{\beta }^{\ast }(g_{L,\eta })},$ respectively, and each of the
above inequalities does not hold in larger disks. By using the inequalities (with $-$ and $+$ in the brackets) in
\eqref{ineqlem} for the case $\lambda=0$ we have that
\begin{align*}
\real\left[\frac{zf_{L,\eta }^{\prime }(z)}{f_{L,\eta }(z)}\right]&=1+
\frac{\eta \real z}{\left( L+1\right) ^{2}}-\frac{1}{L+1}
\sum\limits_{n\geq 1}\real\left[\frac{z^{2}}{x_{L,\eta ,n}(x_{L,\eta ,n}-z)}+\frac{z^{2}}{y_{L,\eta ,n}(y_{L,\eta ,n}-z)}\right] \\
&=1+
\frac{\eta \real z}{\left( L+1\right) ^{2}}-\frac{1}{L+1}
\sum\limits_{n\geq 1}\real\left[\frac{z^{2}}{x_{L,\eta ,n}(x_{L,\eta ,n}-z)}+\frac{z^{2}}{\left|y_{L,\eta ,n}\right|(\left|y_{L,\eta ,n}\right|+z)}\right] \\
&\geq 1+\frac{\eta \left\vert z\right\vert }{\left( L+1\right) ^{2}}-\frac{1%
}{L+1}\sum\limits_{n\geq 1}\left[\frac{\left\vert z\right\vert ^{2}}{x_{L,\eta ,n}(x_{L,\eta ,n}-\left\vert z\right\vert )}+\frac{\left\vert z\right\vert ^{2}}{\left|y_{L,\eta ,n}\right|(\left|y_{L,\eta ,n}\right|+\left\vert z\right\vert )}\right]\\
&=1+\frac{\eta \left\vert z\right\vert }{\left( L+1\right) ^{2}}-\frac{1%
}{L+1}\sum\limits_{n\geq 1}\left[\frac{\left\vert z\right\vert ^{2}}{x_{L,\eta ,n}(x_{L,\eta ,n}-\left\vert z\right\vert )}+\frac{\left\vert z\right\vert ^{2}}{y_{L,\eta ,n}(y_{L,\eta ,n}-\left\vert z\right\vert )}\right]\\
&=\frac{\left\vert z\right\vert f_{L,\eta }^{\prime }(\left\vert
z\right\vert )}{f_{L,\eta }(\left\vert z\right\vert )}
\end{align*}
and
\begin{align*}
\real\left[\frac{zg_{L,\eta }^{\prime }(z)}{g_{L,\eta }(z)}\right]&=1+
\frac{\eta \real z}{L+1}-\sum\limits_{n\geq 1}\real\left[\frac{z^{2}}{x_{L,\eta ,n}(x_{L,\eta ,n}-z)}+\frac{z^{2}}{y_{L,\eta ,n}(y_{L,\eta ,n}-z)}\right] \\
&=1+
\frac{\eta \real z}{L+1}-\sum\limits_{n\geq 1}\real\left[\frac{z^{2}}{x_{L,\eta ,n}(x_{L,\eta ,n}-z)}+\frac{z^{2}}{\left|y_{L,\eta ,n}\right|(\left|y_{L,\eta ,n}\right|+z)}\right] \\
&\geq 1+\frac{\eta \left\vert z\right\vert }{L+1}-\sum\limits_{n\geq 1}\left[\frac{\left\vert z\right\vert ^{2}}{x_{L,\eta ,n}(x_{L,\eta ,n}-\left\vert z\right\vert )}+\frac{\left\vert z\right\vert ^{2}}{\left|y_{L,\eta ,n}\right|(\left|y_{L,\eta ,n}\right|+\left\vert z\right\vert )}\right]\\
&=1+\frac{\eta \left\vert z\right\vert }{L+1}-\sum\limits_{n\geq 1}\left[\frac{\left\vert z\right\vert ^{2}}{x_{L,\eta ,n}(x_{L,\eta ,n}-\left\vert z\right\vert )}+\frac{\left\vert z\right\vert ^{2}}{y_{L,\eta ,n}(y_{L,\eta ,n}-\left\vert z\right\vert )}\right]\\
&=\frac{\left\vert z\right\vert g_{L,\eta }^{\prime }(\left\vert
z\right\vert )}{g_{L,\eta }(\left\vert z\right\vert )},
\end{align*}
where $|z|<\omega=\min\left\{x_{L,\eta,1},\left|y_{L,\eta,1}\right|\right\},$ and equalities are attained only when $z=\left\vert z\right\vert =r.$ The
latter inequalities and minimum principle for harmonic functions imply that
the corresponding inequalities in (\ref{eq6}) hold if and only if $
\left\vert z\right\vert <\alpha _{L,\eta }$ and $\left\vert
z\right\vert <\delta _{L,\eta },$ respectively, where $\alpha _{L,\eta }$
and $\delta _{L,\eta }$ are the smallest positive roots of the equations%
\begin{equation}\label{starfg}
\frac{rf_{L,\eta }^{\prime }(r)}{f_{L,\eta }(r)}=\beta \text{ \ \ and \ \ }%
\frac{rg_{L,\eta }^{\prime }(r)}{g_{L,\eta }(r)}=\beta
\end{equation}
which are equivalent to
\begin{equation*}
rF_{L,\eta }^{\prime }(r)-\beta (L+1)F_{L,\eta }(r)=0\text{ \ \ and \ \ }%
(L+\beta )F_{L,\eta }(r)-rF_{L,\eta }^{\prime }(r)=0.
\end{equation*}
In other words, we proved that
$$\inf_{z\in\mathbb{D}_r}\left[\frac{zf_{L,\eta }^{\prime }(z)}{f_{L,\eta }(z)}\right]=\frac{rf_{L,\eta }^{\prime }(r)}{f_{L,\eta }(r)}=F(r)\quad \mbox{and}
\quad \inf_{z\in\mathbb{D}_r}\left[\frac{zg_{L,\eta }^{\prime }(z)}{g_{L,\eta }(z)}\right]=\frac{rg_{L,\eta }^{\prime }(r)}{g_{L,\eta }(r)}=G(r),$$
and since the real functions $F,G:(0,\omega)\to\mathbb{R}$ are decreasing, and take the limits
$$\lim_{r\searrow 0}F(r)=\lim_{r\searrow 0}G(r)=1\quad \mbox{and}\quad \lim_{r\nearrow x_{L,\eta,1}}F(r)=\lim_{r\searrow y_{L,\eta,1}}F(r)=
\lim_{r\nearrow x_{L,\eta,1}}G(r)=\lim_{r\searrow y_{L,\eta,1}}G(r)=-\infty,$$
it follows that the inequalities in \eqref{eq6} indeed hold for $z\in\mathbb{D}_{\alpha_{L,\eta}}$ and $z\in\mathbb{D}_{\delta_{L,\eta}},$ respectively if and only if the real numbers $\alpha_{L,\eta}$ and $\delta_{L,\eta}$ are the smallest roots of the corresponding equations in \eqref{starfg}.
\end{proof}

\begin{proof}[\bf Proof of Theorem \ref{thuniv}]
Since the coefficients of the normalized regular Coulomb wave functions $f_{L,\eta}$ and $g_{L,\eta}$ are real, according to Wilf \cite{wilf} it follows that their radii of convexity are less or equal than their radii of starlikeness, which are less or equal than their radii of univalence. Thus, we just need to
show that the radii of univalence are less or equal than the corresponding radii of starlikeness. On the other hand, following the proof of Theorem \ref{t2} it is clear that for $L>-1,$ $\eta\leq0$ and $|z|<r$ we have
$$\real\left[\frac{zf_{L,\eta}'(z)}{f_{L,\eta}(z)}\right]>\frac{rf_{L,\eta}'(r)}{f_{L,\eta}(r)}\quad \mbox{and}\quad \real\left[\frac{zg_{L,\eta}'(z)}{g_{L,\eta}(z)}\right]>\frac{rg_{L,\eta}'(r)}{g_{L,\eta}(r)},$$
and each quantity on the right-hand sides of the above inequalities remain positive until the first positive zero of $f_{L,\eta}'$ and $g_{L,\eta}'$ is reached, that is, as long as $r\leq r^{\ast}(f_{L,\eta})$ and $r\leq r^{\ast}(g_{L,\eta}),$ respectively. This shows that indeed the radii of univalence and starlikeness coincide.

Now, recall that the infinite sum of the regular Coulomb wave
function $F_{L,\eta}$ is
$$F_{L,\eta }(z)=C_{L}(\eta )\sum_{n\geq0}a_{L,n}z^{n+L+1}.$$
Thus, we obtain
\begin{equation*}
F_{L,\eta }^{\prime }(z)=C_{L}(\eta )z^{L}\sum_{n\geq0}(n+L+1)a_{L,n}z^{n}
\end{equation*}%
and since the growth order $\tau_{C}$ of the entire regular Coulomb wave function
$z\mapsto C_{L}^{-1}(\eta )z^{-L-1}F_{L,\eta}(z)$
satisfies (see \cite{Ba4})
\begin{equation}
1\leq \tau _{C}=\underset{n\rightarrow \infty }{\lim \sup }\frac{n\log n}{%
-\log |a_{L,n}|}<2,  \label{eq7}
\end{equation}
it follows that
\begin{equation*}
\underset{n\rightarrow \infty }{\lim \sup }\frac{n\log n}{-\log \left[(n+L+1)|a_{L,n}|\right]}=\underset{n\rightarrow \infty }{\lim \sup }%
\frac{n\log n}{-\log (n+L+1)-\log |a_{L,n}|}=\underset{n\rightarrow
\infty }{\lim \sup }\frac{n\log n}{-\log |a_{L,n}|}.
\end{equation*}%
Thus, in view of (\ref{eq7}), and by using the Hadamard theorem concerning the
canonical Weierstrassian representation of entire functions we can write the
infinite product representation
of $F_{L,\eta }^{\prime }$ as follows
\begin{equation}
F_{L,\eta }^{\prime }(z)=(L+1)C_{L}(\eta )z^{L}e^{\frac{\eta(L+2) z}{(L+1)^2}%
}\prod\limits_{n\geq 1}\left( 1-\frac{z}{\rho _{L,\eta ,n}^{\prime }}%
\right) e^{\frac{z}{\rho _{L,\eta ,n}^{\prime }}},  \label{eq8}
\end{equation}
where $\rho _{L,\eta ,n}^{\prime }$ is the $n$th zero of the function $%
F_{L,\eta }^{\prime }(z).$ In other words, the derivative of the Coulomb wave function has the same
growth order as the function itself (this follows also from the theory of entire functions) and consequently
it has a similar infinite product representation. Summarizing, the expression
$$\Theta_{L,\eta}(z)=(L+1)^{-1}C_L(\eta)z^{-L}F_{L,\eta}'(z)$$
can be written as
$$\Theta_{L,\eta}(z)=e^{\frac{\eta(L+2) z}{(L+1)^2}%
}\prod\limits_{n\geq 1}\left( 1-\frac{z}{\rho _{L,\eta ,n}^{\prime }}%
\right) e^{\frac{z}{\rho _{L,\eta ,n}^{\prime }}}=\sum_{n\geq0}\frac{n+L+1}{L+1}a_{L,n}z^n,$$
and taking its logarithmic derivative we obtain that
$$\frac{\Theta'_{L,\eta}(z)}{\Theta_{L,\eta}(z)}=\frac{\eta(L+2)}{(L+1)^2}+\sigma_1-\sum_{n\geq1}\sum_{m\geq0}\frac{z^m}{({\rho_{L,\eta,n}'})^{m+1}}=
\frac{\eta(L+2)}{(L+1)^2}-\sum_{m\geq1}\sigma_{m+1}z^m$$
and
$$\frac{\Theta'_{L,\eta}(z)}{\Theta_{L,\eta}(z)}=\frac{\sum\limits_{n\geq0}\frac{(n+1)(n+L+2)}{L+1}a_{L,n+1}z^n}{\sum\limits_{n\geq0}\frac{n+L+1}{L+1}a_{L,n}z^n},$$
where $|z|<|\rho_{L,\eta,1}'|$ and $\sigma_m=\sum_{n\geq1}{(\rho_{L,\eta,n}')}^{-m}$ is the Rayleigh sum of the zeros of the derivative of the Coulomb wave function.
By equating the coefficients of the power series we can find the Rayleigh sums $\sigma_m$ term by term and by using the Euler-Rayleigh inequalities
$$\sigma_m^{-\frac{1}{m}}<|\rho_{L,\eta,1}'|<\frac{\sigma_m}{\sigma_{m+1}}$$
we can found lower and upper bounds for the radius of starlikeness (and univalence) of the function $f_{L,\eta}.$ To prove the inequalities \eqref{ray1} and \eqref{ray2} we shall use the above Euler-Rayleigh inequalities for $m=2.$ However, we remark here that in the Euler-Rayleigh inequalities the lower bounds increase and upper bounds decrease to $\rho_{L,\eta,1}'$ as $m\to\infty.$ Thus, it would be possible to improve the inequalities \eqref{ray1} and \eqref{ray2} by using higher order Euler-Rayleigh inequalities. We note that
$$\sigma_2=-\lim_{z\to0}\left[\frac{\Theta'_{L,\eta}(z)}{\Theta_{L,\eta}(z)}\right]'=\left(\frac{L+2}{L+1}a_{L,1}\right)^2-2\frac{L+3}{L+1}a_{L,2}$$
and
$$\sigma_3=-\lim_{z\to0}\left[\frac{\Theta'_{L,\eta}(z)}{\Theta_{L,\eta}(z)}\right]''=
\frac{\eta(L+2)}{(L+1)^2}\frac{L+3}{L+1}a_{L,2}-\sigma_2\frac{L+2}{L+1}a_{L,1}-\frac{3(L+4)}{L+1}a_{L,3},$$
which can be rewritten as
$$\sigma_2=\frac{L^4+6L^3+(\eta^2+12)L^2+2(3\eta^2+5)L+3(2\eta^2+1)}{(L+1)^4(2L+3)}$$
and
$$\sigma_3=\frac{\eta(4L^7+\eta_1L^6+\eta_2L^5+8\eta_3L^4+2\eta_4L^3+\eta_5 L^2+9(5-28\eta^2)L-72\eta^2)}{2(L+1)^6(L+2)(2L+3)}.$$

Now, we focus on the radius of starlikeness of the normalized regular Coulomb wave function $g_{L,\eta},$ which by means of \eqref{prodco} can be written as
$$g_{L,\eta}(z)=ze^{\frac{\eta z}{L+1}}\prod_{n\geq 1}\left(1-\frac{z}{\rho_{L,\eta,n}}\right)e^{\frac{z}{\rho_{L,\eta,n}}}.$$
Since for $L>-1$ and $\eta\in\mathbb{R}$ the zeros of the regular Coulomb wave functions are real it follows that the function $g_{L,\eta}$ belongs to the Laguerre-P\'olya class $\mathcal{LP}$ of real entire functions, which are uniform limits of real polynomials whose all zeros are real. Now, since the Laguerre-P\'olya class $\mathcal{LP}$ is closed under differentiation, it follows that $g_{L,\eta}'$ belongs also to the Laguerre-P\'olya class and hence all of its zeros are real. On the other hand, by using the fact the growth order of an entire function and its derivative coincide, in view of
Hadamard theorem we immediately have that
\begin{equation}\label{prodgder}g_{L,\eta}'(z)=e^{\frac{2\eta z}{L+1}}\prod_{n\geq 1}\left(1-\frac{z}{\xi_{L,\eta,n}}\right)e^{\frac{z}{\xi_{L,\eta,n}}},\end{equation}
where $\xi_{L,\eta,n}$ stands for the $n$th zero of the function $g_{L,\eta}',$ that is of $r\mapsto rF_{L,\eta}'(r)-LF_{L,\eta}(r).$ Thus, by using this infinite product we obtain
$$\frac{g_{L,\eta}''(z)}{g_{L,\eta}'(z)}=\frac{2\eta}{L+1}+\varsigma_1-\sum_{n\geq1}\sum_{k\geq0}\frac{z^k}{\xi_{L,\eta,n}^{k+1}}=
\frac{2\eta}{L+1}-\sum_{k\geq1}\varsigma_{k+1}z^k,$$
where $|z|<|\xi_{L,\eta,1}|$ and $\varsigma_k=\sum_{n\geq1}\xi_{L,\eta,n}^{-k}$ is the corresponding Rayleigh sum. Equating the corresponding coefficients we obtain that
$$\varsigma_2=\frac{4\eta}{L+1}a_{L,1}-6a_{L,2}\quad \mbox{and}\quad \varsigma_3=\frac{6\eta}{L+1}a_{L,2}-2\varsigma_2a_{L,1}-12a_{L,3},$$
which can be written as
$$\varsigma_2=\frac{3L^2+2(\eta^2+3)L+3(2\eta^2+1)}{(L+1)^2(2L+3)}$$
and
$$\varsigma_3=\frac{\eta(8L^4+(31-16\eta^2)L^3+2(19-26\eta^2)L^2+3(5-22\eta^2)L-36\eta^2)}{(L+1)^3(L+2)(2L+3)}.$$
Thus, by using the Euler-Rayleigh inequalities
$$\varsigma_m^{-\frac{1}{m}}<|\xi_{L,\eta,1}|<\frac{\varsigma_m}{\varsigma_{m+1}}$$
for the particular case $m=2,$ we obtain the estimates in \eqref{ray3} and \eqref{ray4}.
\end{proof}

\begin{proof}[\bf Proof of Theorem \ref{thconv}]
We note that
\begin{equation*}
1+\frac{zf_{L,\eta }^{\prime \prime }(z)}{f_{L,\eta }^{\prime }(z)}=1+\frac{%
zF_{L,\eta }^{\prime \prime }(z)}{F_{L,\eta }^{\prime }(z)}-\frac{L}{L+1}\frac{zF_{L,\eta }^{\prime }(z)}{F_{L,\eta }(z)},
\end{equation*}
and since by means of (\ref{eq8}) we have
$$
1+\frac{zF_{L,\eta }^{\prime \prime }(z)}{F_{L,\eta }^{\prime }(z)}=1+L+%
\frac{\eta(L+2) z}{(L+1)^2}-\sum\limits_{n\geq 1}\frac{z^{2}}{\rho _{L,\eta
,n}^{\prime }(\rho _{L,\eta ,n}^{\prime }-z)},  \label{eq9}
$$
it follows that
\begin{equation*}
1+\frac{zf_{L,\eta }^{\prime \prime }(z)}{f_{L,\eta }^{\prime }(z)}=1+\frac{%
2\eta z}{\left( L+1\right) ^{2}}+\frac{L}{L+1}
\sum\limits_{n\geq 1}\frac{z^{2}}{\rho _{L,\eta ,n}(\rho _{L,\eta ,n}-z)}%
-\sum\limits_{n\geq 1}\frac{z^{2}}{\rho _{L,\eta ,n}^{\prime }(\rho
_{L,\eta ,n}^{\prime }-z)}.
\end{equation*}
On the other hand, we know that (see \cite{St}) if $L>-1$ and $\eta \in\mathbb{R},$ then zeros
of the regular Coulomb wave function $F_{L,\eta}$ are real. Moreover, if
$L>-\frac{1}{2}$ and $\eta \in\mathbb{R}$, then the zeros of $F_{L,\eta }$ and
$F_{L,\eta }^{\prime }$ are interlacing (see \cite{Ba4} for more details). Denoting by $x_{L,\eta,n}'$ and $y_{L,\eta,n}'$
the $n$th positive and negative zeros of the function $F_{L,\eta}'$ we obtain for $L\in\left(-\frac{1}{2},0\right]$ and $|z|<\omega'=\min\{x_{L,\eta,1}',y_{L,\eta,1}'\}$
\begin{align*}
\real\left[1+\frac{zf_{L,\eta}''(z)}{f_{L,\eta }'(z)}\right]&=1+
\frac{2\eta \real z}{\left( L+1\right) ^{2}}+\frac{L}{L+1}
\sum\limits_{n\geq 1}\real\left[\frac{z^{2}}{x_{L,\eta ,n}(x_{L,\eta ,n}-z)}+\frac{z^{2}}{y_{L,\eta ,n}(y_{L,\eta ,n}-z)}\right]\\&-
\sum\limits_{n\geq 1}\real\left[\frac{z^{2}}{x_{L,\eta ,n}'(x_{L,\eta ,n}'-z)}+\frac{z^{2}}{y_{L,\eta ,n}'(y_{L,\eta ,n}'-z)}\right]
\end{align*}
and proceeding exactly as in the proof of Theorem \ref{t2}, we obtain that
$$\real\left[1+\frac{zf_{L,\eta}''(z)}{f_{L,\eta }'(z)}\right]\geq 1+\frac{|z|f_{L,\eta}''(|z|)}{f_{L,\eta }'(|z|)}$$
for all $L\in\left(-\frac{1}{2},0\right]$ and $|z|<\omega'.$ Moreover, observe that if we use the
inequality \eqref{ineqlem} we get that the above
inequality is also valid when $L>0.$ Here we used that the
zeros $\rho _{L,\eta ,n}$ and $\rho _{L,\eta ,n}^{\prime }$ interlace for $%
L>-\frac{1}{2},$ according to \cite[Theorem 5]{Ba4}. Thus, for $r\in (0,\omega')$ we have
$$\underset{z\in \mathbb{D}_r}{\inf }\left\{ \real\left[ 1+%
\frac{zf_{L,\eta }^{\prime \prime }(z)}{f_{L,\eta }^{\prime }(z)}\right]
\right\} =1+\frac{rf_{L,\eta }^{\prime \prime }(r)}{f_{L,\eta }^{\prime }(r)}.
$$
On the other hand, the function $\mathbb{F}_{L,\eta }:$ $\left(0,\omega'\right)\to \mathbb{R},$ defined by
\begin{equation*}
\mathbb{F}_{L,\eta }(r)=1+\frac{rf_{L,\eta }^{\prime \prime }(r)}{f_{L,\eta
}^{\prime }(r)}
\end{equation*}%
is strictly decreasing for all $\eta\leq 0$ and $L>-\frac{1}{2}.$ Namely, we have
\begin{equation*}
\mathbb{F}_{L,\eta }^{\prime }(r)=\frac{2\eta }{\left( L+1\right) ^{2}}+\frac{L}{L+1} \sum\limits_{n\geq 1}\frac{r(2\rho _{L,\eta
,n}-r)}{\rho _{L,\eta ,n}(\rho _{L,\eta ,n}-r)^{2}}-\sum\limits_{n\geq 1}%
\frac{r(2\rho _{L,\eta ,n}^{\prime }-r)}{\rho _{L,\eta ,n}^{\prime }(\rho
_{L,\eta ,n}^{\prime }-r)^{2}},
\end{equation*}
and this is clearly negative for $\eta\leq0,$ $L\in\left(-\frac{1}{2},0\right]$ and $r<\omega'.$ Moreover, when $\eta\leq0$ and $L>0$ we have
\begin{equation*}
\mathbb{F}_{L,\eta }^{\prime }(r)<\sum\limits_{n\geq 1}\frac{r(2\rho _{L,\eta
,n}-r)}{\rho _{L,\eta ,n}(\rho _{L,\eta ,n}-r)^{2}}-\sum\limits_{n\geq 1}%
\frac{r(2\rho _{L,\eta ,n}^{\prime }-r)}{\rho _{L,\eta ,n}^{\prime }(\rho
_{L,\eta ,n}^{\prime }-r)^{2}}<0,
\end{equation*}
where we used again the fact that the
zeros $\rho _{L,\eta ,n}$ and $\rho _{L,\eta ,n}^{\prime }$ interlace for $%
L>-\frac{1}{2},$ and the functions $a\mapsto \displaystyle\frac{2a\pm r}{a(a\pm r)^2}$ are decreasing. On the other hand, we have the limits
$\lim_{r\searrow0}\mathbb{F}_{L,\eta }(r)=1$ and $\lim_{r\nearrow\omega'}\mathbb{F}_{L,\eta }(r)=-\infty.$ Consequently, the equation $\mathbb{F}_{L,\eta }(r)=\beta$ has a unique root and this equation is equivalent to \eqref{eqconv}.

Now, we focus on the radius of convexity of the normalized regular Coulomb wave function $g_{L,\eta}.$ Since the infinite product in \eqref{prodgder} is absolutely convergent on compact subsets of the complex
plane, we can take the logarithmic derivative of both sides of the equation \eqref{prodgder} to obtain
$$1+\frac{zg_{L,\eta}''(z)}{g_{L,\eta}'(z)}=1+\frac{2\eta z}{L+1}-\sum_{n\geq 1}\frac{z^2}{\xi_{L,\eta,n}(\xi_{L,\eta,n}-z)}.$$
By using similar steps as in the proof of Theorem \ref{t2} we conclude that
$$\real\left[1+\frac{zg_{L,\eta}''(z)}{g_{L,\eta}'(z)}\right]\geq 1+\frac{|z|g_{L,\eta}''(|z|)}{g_{L,\eta}'(|z|)}$$
where $|z|<\varepsilon=\min\left\{\alpha_{L,\eta,1},\left|\beta_{L,\eta,1}\right|\right\},$ and equality is attained only when $z=\left\vert z\right\vert =r.$ Here $\alpha_{L,\eta,n}$ and $\beta_{L,\eta,n}$ denote the $n$th positive and negative zero of $g_{L,\eta}'.$ The above inequality and minimum principle for harmonic functions imply that the inequality
$$\real\left[1+\frac{zg_{L,\eta}''(z)}{g_{L,\eta}'(z)}\right]>\beta$$
holds if and only if $\left\vert z\right\vert <\varsigma_{L,\eta },$ where $\varsigma_{L,\eta }$
is the smallest positive root of the equations
$$1+\frac{rg_{L,\eta}''(r)}{g_{L,\eta}'(r)}=\beta.$$
Here we used the fact that the function
$$r\mapsto 1+\frac{rg_{L,\eta}''(r)}{g_{L,\eta}'(r)}=1+\frac{2\eta r}{L+1}-\sum_{n\geq 1}\frac{r^2}{\xi_{L,\eta,n}(\xi_{L,\eta,n}-r)}$$
maps $(0,\varepsilon)$ into $(-\infty,1)$ and it is decreasing.
\end{proof}

\begin{proof}[\bf Proof of Theorem \ref{thsub}]
Since the regular Coulomb wave function $F_{L,\eta}$ satisfies the Coulomb differential equation
$$z^2w''(z)+\left[z^2-2\eta z-L(L+1)\right]w(z)=0$$
it follows that $g_{L,\eta}$ satisfies the next homogeneous second order linear differential equation
$$z^2w''(z)+2Lzw'(z)+\left( z^{2}-2\eta
z-2L\right)w(z)=0$$
and consequently
\begin{equation}
z^{2}g_{L,\eta}^{\prime \prime }(z)+2Lzg_{L,\eta}^{\prime }(z))+\left( z^{2}-2\eta
z-2L\right)g_{L,\eta}(z)=0,  \label{equa1}
\end{equation}%
where $z\in \mathbb{D},$ $\eta ,L\in\mathbb{C}.$ If we consider the expression $\Psi (r,s,t;z)=t+2Ls+\left( z^{2}-2\eta z-2L\right) r$ and $%
\mathbb{C} \supseteq \Omega=\{0\}$, then the equation (\ref{equa1}) implies $\Psi
(g_{L,\eta}(z),zg_{L,\eta}^{^{\prime }}(z),z^{2}g_{L,\eta}^{^{\prime \prime }}(z);z)\in \Omega$
for all $z\in \mathbb{D}.$ Now, we shall use Lemma \ref{l1} to prove that $\real
g_{L,\eta}(z)>0$ for all $z\in \mathbb{D}.$ If we put $z=x+iy,$ where $x,y\in(-1,1)
,$ then
\begin{equation*}
\real\Psi \left( \rho i,\sigma ,\mu +i\nu;x+iy\right) =\left( \mu +\sigma
\right) +\left( 2\real L-1\right) \sigma -\left( 2xy-2y\real\eta -2x%
\imag\eta -2\imag L\right) \rho
\end{equation*}%
for all $\rho ,\sigma ,\mu ,\nu\in
\mathbb{R}.$ Let $\rho ,\sigma ,\mu ,\nu\in
\mathbb{R}$ and $L\in\mathbb{C}$ satisfy $\mu +\sigma \leq 0,$ $\sigma \leq -(1+\rho ^{2})/2$ and $\real L\geq \frac{1}{2}.$ We have%
\begin{equation*}
\real\Psi \left( \rho i,\sigma ,\mu +i\nu;x+iy\right) \leq -\left(\real L-\frac{1}{2}\right)\rho ^{2}-2\left(xy-y\real\eta-x\imag
\eta-\imag L\right)\rho -\left(\real L-\frac{1}{2}\right)=Q_1(\rho
).
\end{equation*}%
The discriminant $\Delta_1$ of the quadratic form $Q_1(\rho )$ is%
\begin{equation*}
\Delta_1 =4\left(y\real\eta+x\imag\eta+\imag L-xy\right)
^{2}-\left( 2\real L-1\right) ^{2}.
\end{equation*}%
By using the well-known Cauchy-Bunyakovsky-Schwarz inequality we have%
\begin{equation*}
y\real\eta +x\imag\eta \leq \left\vert y\real\eta +x\imag
\eta \right\vert \leq \sqrt{\left( \real\eta \right) ^{2}+\left( \imag\eta \right) ^{2}}\sqrt{x^{2}+y^{2}}<\left\vert \eta \right\vert ,
\end{equation*}%
and since $\imag L-xy>\imag L-1\geq 0$ we get that
\begin{align*}
\Delta_1 &=4\left[\left(y\real\eta +x\imag\eta\right)^2+2\left(y\real\eta +x\imag\eta\right)\left(\imag L-xy\right)+\left(\imag L-xy\right)^2\right]-\left( 2\real L-1\right) ^{2}\\
&\leq 4\left[|\eta|^2+2|\eta|(\imag L+1)+\left(\imag L+1\right)^2\right]-\left( 2\real L-1\right) ^{2}\\
&\leq 4\left(1+\imag L+\left\vert \eta \right\vert \right)
^{2}-\left( 2\real L-1\right) ^{2}\leq 0.
\end{align*}
Thus, the quadratic form $Q_1(\rho )$ is strictly negative and consequently we
have $$\real\Psi \left( \rho i,\sigma ,\mu +i\nu;x+iy\right) <0.$$ By Lemma %
\ref{l1} we conclude that $\real g_{L,\eta}(z)>0$ for all $z\in \mathbb{D}.$

Now, define the function $q_{L,\eta}:\mathbb{D}\to\mathbb{C}$ by
\begin{equation*}
q_{L,\eta}(z)=\frac{zg_{L,\eta}^{\prime}(z)}{g_{L,\eta}(z)}.
\end{equation*}%
The function $q_{L,\eta}$ is analytic in $\mathbb{D}$ and $q_{L,\eta}(0)=1.$ Suppose that $z\neq 0.$ We
know that $g_{L,\eta}(z)\neq 0$ and therefore if we divide both sides of equation (%
\ref{equa1}) with $g_{L,\eta}(z),$ we obtain%
\begin{equation*}
\frac{z^{2}g_{L,\eta}^{\prime \prime }(z)}{g_{L,\eta}(z)}+2L\frac{zg_{L,\eta}^{\prime }(z)}{%
g_{L,\eta}(z)}+\left(z^{2}-2\eta z-2L\right) =0
\end{equation*}%
or%
\begin{equation}
\left[ \frac{zg_{L,\eta}^{\prime \prime }(z)}{g_{L,\eta}^{\prime }(z)}\right] \left[
\frac{zg_{L,\eta}^{\prime }(z)}{g_{L,\eta}(z)}\right] +2L\left[ \frac{zg_{L,\eta}^{\prime
}(z)}{g_{L,\eta}(z)}\right] +\left( z^{2}-2\eta z-2L\right) =0.  \label{eq2}
\end{equation}%
Differentiating logarithmically the expression $q_{L,\eta}(z),$ we obtain%
\begin{equation*}
\frac{zg_{L,\eta}^{\prime \prime }(z)}{g_{L,\eta}^{\prime }(z)}=\frac{zq^{\prime
}(z)+q^{2}(z)-q(z)}{q(z)}.
\end{equation*}%
In view of (\ref{eq2}) this result reveals that $q_{L,\eta}$ satisfies the following
differential equation:%
\begin{equation}
zq_{L,\eta}^{\prime }(z)+q_{L,\eta}^{2}(z)+(2L-1)q_{L,\eta}(z)+\left( z^{2}-2\eta z-2L\right) =0.
\label{eq3}
\end{equation}%
If we use $\Psi (r,s;z)=s+r^{2}+(2L-1)r+\left( z^{2}-2\eta z-2L\right) $ and
$\Omega=\{0\},$ then (\ref{eq3}) implies that we have $\Psi (q_{L,\eta}(z),zq_{L,\eta}^{\prime }(z);z)\in \Omega$ for
all $z\in \mathbb{D}.$ Now we use Lemma \ref{l1} to prove that $\real q_{L,\eta}(z)>0$ for
all $z\in \mathbb{D}.$ For $z=x+iy\in \mathbb{D}$ and $\rho ,\sigma \in\mathbb{R}$ satisfying $\sigma \leq -(1+\rho ^{2})/2$ we
obtain
\begin{align*}
\real\Psi \left( \rho i,\sigma ;x+iy\right)  &=\sigma -\rho
^{2}+x^{2}-y^{2}+2\rho\imag L-2x\real\eta +2y\imag\eta -2\real L \\
&\leq -\frac{3}{2}\rho ^{2}+2\rho \imag L-2x\real\eta +2y\imag%
\eta -2\real L+\frac{1}{2}=Q_2(\rho)
\end{align*}
The discriminant $\Delta_2$ of the quadratic form $Q_2(\rho )$ is%
\begin{equation*}
\Delta_2 =4(\imag L)^{2}+12\left(y\imag \eta -x\real\eta -\real L+\frac{1}{4}\right).
\end{equation*}%
By the Cauchy-Bunyakovsky-Schwarz inequality we have%
\begin{equation*}
y\imag\eta -x\real\eta \leq \left\vert y\imag\eta -x\real
\eta \right\vert \leq \sqrt{\left(\real\eta \right) ^{2}+\left(\imag\eta \right) ^{2}}\sqrt{x^{2}+y^{2}}<\left\vert \eta \right\vert .
\end{equation*}%
Therefore we have%
\begin{equation*}
\frac{\Delta_2}{4}\leq \left(\imag L\right) ^{2}+3\left(\left\vert \eta
\right\vert -\real L+\frac{1}{4}\right) \leq 0.
\end{equation*}%
Thus, the quadratic form $Q_2(\rho)$ is strictly negative and consequently we
have $\real\Psi \left(\rho i,\sigma ;x+iy\right) <0.$ By Lemma \ref{l1}
we conclude that $\real q_{L,\eta}(z)>0$ for all $z\in \mathbb{D},$ which show that $g_{L,\eta}$
is starlike in $\mathbb{D}.$
\end{proof}

\end{document}